\documentclass[12pt,reqno,a4paper]{amsart}
\usepackage{extsizes}

\usepackage{amsthm}
\usepackage{amstext}
\usepackage{amssymb}
\usepackage{mathrsfs}
\usepackage{mathabx}
\usepackage[english]{babel}
\usepackage{enumerate}
\usepackage{xfrac}
\usepackage{mathabx}
\usepackage{mathtools}
\usepackage{xcolor}
\usepackage[all]{xy}

\usepackage[T1]{fontenc}

\theoremstyle{plain}
\newtheorem{theorem}{Theorem}[section]

\newtheorem{proposition}[theorem]{Proposition}

\newtheorem{lemma}[theorem]{Lemma}
\newtheorem{propandef}[theorem]{Proposition and Definition}

\theoremstyle{definition}

\newtheorem{remark}[theorem]{Remark}

\theoremstyle{remark}

\numberwithin{equation}{section}

\newcommand{\Id}{\, d}

\newcommand{\IR}{\mathbb{R}}

\def\bf{\mathbf}

\newcommand{\loc}{\mathrm{loc}}

\newcommand{\IMM}{\mathscr{B}}

\newcommand{\IN}{\mathbb{N}}

\newcommand\newdot{{\kern.8pt\cdot\kern.8pt}}
\def\nbull{{\raise1.5pt\hbox{\bf .}}}

\title[]{Locally convex aspects of the Kato and the Dynkin class on manifolds}%

\author{Batu G\"uneysu}

\author{Kazuhiro Kuwae}

\thanks{Fakultät für Mathematik, Technische Universität Chemnitz,
	09126 Chemnitz, Germany, {\tt batu.gueneysu@math} {\tt tu-chemnitz.de}\\
Department of Applied Mathematics, Fukuoka University,
Fukuoka 814-0180, Japan, {\tt kuwae@fukuoka-u.} {\tt ac.jp}
}

\begin{document}

\maketitle

\begin{abstract} Given a Riemannian manifold $(X,g)$ we consider the Kato class $\mathcal{K}(X,g)$ and the Dynkin class $\widetilde{\mathcal{K}}_\loc(X,g)$, and well as their local counterparts as Fr\'{e}chet spaces. Based on recent results by Carron, Mondello and Tewodrose, we show that if $\dim(X)\geq 2$ and if spectral negative part $\sigma^-_g$ of the Ricci curvature is in $L^q_{\phi_g}(X,g)+L^\infty(X,g)$ for some $q> 1$, the assumption $\sigma^-_g\in \mathcal{K}(X,g)$ is equivalent to a Gaussian upper bound of the heat kernel for small times together with uniform local volume doubling (if in both cases $q$ is large enough). Here $L^q_{\phi_g}(X,g)$ is the $L^q$-space which is weighted with the inverse volume function.\\
By establishing a localization result for the Dynkin norm, we prove that the local Kato class and the local Dynkin class do not depend on the chosen Riemannian metric and thus can be defined as Fr\'{e}chet spaces on arbitrary smooth manifolds. Moreover, we prove that smooth compactly supported functions are dense in the local Kato class and we use this result to prove that Schr\"odinger semigroups with Kato decomposable potentials are space-time continuous.
\end{abstract}

\section{Introduction}

Kato functions have been introduced in \cite{kato} in Euclidean space in the context of essential self-adjointess of Schr\"odinger operators as follows: one says that a Borel function $w:\IR^m\to \IR$ is in the Kato class $\mathcal{K}(\IR^m)$, if 
\begin{align}
	&w\in L^1_{\mathrm{unif},\loc}(\IR), \quad\text{in case $m=1$},\label{ursp1}\\
        &\lim_{r\to0}\sup_{x\in \IR^m}\int_{|x-y|<r}  \log_+(|x-y|^{-1})|w(y)| dy=0, \quad\text{in case $m=2$},\label{ursp2}\\
        &\lim_{r\to0}\sup_{x\in \IR^m}\int_{|x-y|<r}  |x-y|^{2-m}|w(y)| dy=0, \quad\text{in case $m\geq 3$}.\label{ursp3}
	\end{align}
It has been shown in \cite{aizen} in that $w$ is in $\mathcal{K}(\IR^m)$, if and only if
\begin{align}\label{gau}
\lim_{t\to 0+} \int^t_0\int_{\IR^m} (4\pi s)^{-m/2}e^{ -\frac{|x-y|^2}{4s}   }|w(y)|\Id y \Id s=0,
\end{align}
and one says that $w$ is in the local Kato class $\mathcal{K}_\loc(\IR^m)$, if $1_Kw\in \mathcal{K}(\IR^m)$ for all compact $K\subset \IR^m$. It is straightforward to check that if $m\geq 2$ and $p\in (m/2,\infty)$ then one has an inclusion
\begin{align}
L^p(\IR^m)+L^\infty(\IR^m) \subset \mathcal{K}(\IR^m).
\end{align}
In the world of Schr\"odinger operators, Kato potentials play a crucial role as a global growth condition, as e.g. for every $w\in L^1_\loc(\IR^m)$ with negative part $w^-$ in $\mathcal{K}(\IR^m)$ the Schr\"odinger operator $H^w:=-\Delta+w$ in $L^2(\IR^m)$ turns out to be semibounded from below. On the other hand, local Kato potentials play a crucial role in local regularity theory of Schr\"odinger operators, as e.g. if in the above situation one additionally assumes $w\in \mathcal{K}_\loc(\IR^m)$, then $(t,x)\mapsto e^{-t H^w}\psi(x)$ is jointly continuous for all $\psi\in L^2(\IR^m)$ \cite{hundert}. The continuity (or rather local boundedness) of $e^{-t H^w}\psi(x)$ in $x$ plays a crucial role in the essential self-adjointness of Schr\"odinger operators (see for example the proof of Theorem B.1.6 in \cite{barry} or \cite{batu5}).\vspace{1mm}

Replacing the Euclidean $\IR^m$ with an arbitrary Riemannian manifold $(X,g)$ of dimension $m$, there are in principle two ways of defining the Kato class $\mathcal{K}(X,g)$: one is by replacing in \eqref{ursp1}--\eqref{ursp3} the Euclidean distance with the geodesic distance $\varrho_g$ and the Lebesgue measure with the volume measure $\mu_g$. The second one is by replacing the Gauss-Weierstrass kernel in \eqref{gau} with the (minimal) heat kernel $p_g(s,x,y)$ on $(X,g)$. While under certain assumptions on the geometry of $(X,g)$ both definitions agree \cite{kuwae1}, the heat kernel definition has turned out to be more useful and is considered the standard one (and is taken in this paper, too). Indeed, with $\sigma_g(x)$ the smallest eigenvalue of the Ricci tensor $\mathrm{Ric}_g(x)$, read as an symmetric element of $\mathrm{End}(T_x X)$, in recent years it turned out that the assumption $\sigma_g^-:=-\min(\sigma_g(x) , 0)\in \mathcal{K}(X,g)$ is very useful in geometry, topology and probability \cite{batu2,braun, carron1,carron2,carron4,carron5,carron6,rose}. \\
In this paper, we take the perspective of considering $\mathcal{K}(X,g)$ as a Fr\'{e}chet space. In fact, given a Borel function $w:X\to\IR$ the central data of this paper are the Dynkin norms
$$
\left\|w\right\|_{g,t}:= \sup_{x\in X}\int^t_0 \int_X p_g(s,x,y)|w(y)| \Id\mu_g (y) \Id s\in [0,\infty],\quad t>0.
$$
We show that the space of $w$'s with finite $\left\|w\right\|_{g,t}$ for some/all $t>0$ forms a Fr\'{e}chet space, the \emph{Dynkin class} or \emph{extended Kato class}  $\widetilde{\mathcal{K}}(X,g)$ of $(X,g)$. The Kato class $\mathcal{K}(X,g)$, given by all $w$'s with $\lim_{t\to 0+}\left\|w\right\|_{g,t}=0$, is then a closed subspace of $\widetilde{\mathcal{K}}(X,g)$, thus a Fr\'{e}chet space in itself. The local counterpart $\widetilde{\mathcal{K}}_\loc(X,g)$ is given by all $w$'s such that $1_Kw\in \widetilde{\mathcal{K}}(X,g)$ for all compact $K\subset X$. This again becomes a Fr\'{e}chet space with the family of norms $\left\|w\right\|_{g,t;K}:=\left\|1_Kw\right\|_{g,t}$, and the local Kato class $\mathcal{K}_\loc(X,g)$ (defined analogously) becomes closed subspace of $\widetilde{\mathcal{K}}_\loc(X,g)$.\\
Our main results are as follows:
\begin{itemize}
	\item We show that the supremum over $X$ appearing in $\left\|w\right\|_{g,t;K}$ in fact localizes on the compact $K$ (more generally, this holds for every closed subset).
    \item Based on this localization, we prove that, as locally convex spaces, neither $\widetilde{\mathcal{K}}_\loc(X,g)$ nor $\mathcal{K}_\loc(X,g)$ depend on the Riemannian metric $g$ on $X$; in particular, one can define the Fr\'{e}chet spaces $\widetilde{\mathcal{K}}_\loc(X)$ and $\mathcal{K}_\loc(X)$ and any (smooth) manifold $X$, and under $m\geq 2$ one has $L^q_\loc(X)\subset \mathcal{K}_\loc(X)$ continuously for all $q\in (m/2,\infty)$.
    \item Let $L^q_{\phi_g}(X,g)$ stand for the $L^q$-space w.r.t. the measure $\phi_g(x)\Id\mu_g(x)$, with $1/\phi_g(x)$ the volume of the ball with radius $1$ around $x$. Assume $m\geq 2$. We show that if $\sigma^g_-\in L^q_{\phi_g}(X,g)+L^\infty(X;g)$ for some $q\in (1,\infty)$, then one has $\sigma_g^-\in \mathcal{K}(X,g)$ and $q>N_g$ (where $N_g\geq m$ can be calculated explicitly), if and only if $p_g(t,x,y)$ satisfies an upper Gaussian estimate for small times, $(X,g)$ is locally uniformly $N$-volume doubling for some $N\geq 2$ with $q>N/2$. These global results make a heavy use of the recently established results from \cite{carron3}.
     
    \item Again using the localization result, we show that smooth compactly supported functions are dense in $\mathcal{K}_\loc(X)$ and we use this result to show that $(t,x)\mapsto e^{-t H^w_g}\psi(x)$ is jointly continuous for the Schr\"odinger operator $H_g^w=-\Delta_g+w$ in $L^2(X,g)$ and all $\psi\in L^2(X,g)$, if $w:X\to\IR$ is a potential with $w^+\in \mathcal{K}_\loc(X)$, $w^-\in \mathcal{K}(X,g)$. 
\end{itemize}

\section{Dynkin norms}

Let $X$ be an $m$-dimensional connected manifold\footnote{We understand all (Riemannian) manifolds to be smooth and without boundary in the sequel.}. In the sequel, we understand our function spacees over $\IR$, and any equivalence class of Borel functions on $X$ will be understood with respect to some continuous Borel measure on $X$, that is, a measure which in each chart has a strictly positive continuous density with respect to the Lebesgue measure on $\IR^m$. Also, the support of such an equivalence class will be understood as the essential support with respect some continuous measure. Both notions as well as the locally convex spaces $L^q_\loc(X)$, $q\in [1,\infty]$, do not depend on the particular choice of a continuous measure. \vspace{1mm}

We set
$$
\IMM(X):=\{\text{equivalence classes of Borel measurable functions on $X$}\}.
$$
Let $g$ be a Riemannian metric on $X$. We denote with $\varrho_g(x,y)$ the geodesic distance, with $B_g(x,r)$ the induced open balls, and with $\mu_g$ the Riemannian volume measure, which is nothing but the $m$-dimensional Hausdorff measure on $(X,\varrho_g)$. One also sets $\mu_g(x,r):=\mu_g(B_g(x,r))$. The symbol $\mathrm{inj}_g(x)$ will denote the injectivity radius of $(X,g)$ in $x$.\vspace{1mm}

The (minimal) heat kernel
$$
p_g:(0,\infty)\times X\times X\longrightarrow (0,\infty)
$$
of $(X,g)$ can be defined (see \cite{gri}) as the uniquely determined smooth function, such that for all $t>0$, $f\in L^2(X,g)$, $\mu_g$-a.e. $x\in X$ one has
$$
e^{-t H_g}f(x)=\int_Xp_g(t,x,y)f(y)\Id\mu_g(y),
$$
with $H_g\geq 0$ the Friedrichs realization of the Laplace-Beltrami operator $-\Delta_g$ in $L^2(X,g)$ and $e^{-t H_g}$ is defined via spectral calculus. Recall that $H_g$ is the self-adjoint nonnegative operator induced by the closure $\mathscr{E}_g$ of the symmetric and nonnegative bilinear form in $L^2(X,g)$ defined by 
$$
C_c^\infty(X)\times C_c^\infty(X)\ni (\psi_1,\psi_2)\longmapsto \int_X (-\Delta_g\psi_1) \cdot\psi_2 \Id\mu_g\in \IR.
$$ 
Note that $\mathscr{E}_g$ is a (strongly local) Dirichlet form and that the strict positivity of $p_g$ follows \cite{gri} from the connectedness of $X$.\vspace{1mm}

Given a Borel set $A\subset X$ and a Borel function $\phi:A\to (0,\infty)$, we define $L^q_\phi(A,g)$ to be the $L^q$-space which is defined in terms of Borel measure $\phi  \Id\mu_g$ on $A$. \vspace{1mm}

Given $t>0$ the corresponding \emph{Dynkin norm} of $w\in \IMM(X)$ is defined by
$$
\left\|w\right\|_{g,t}:= \sup_{x\in X}\int^t_0 \int_X p_g(s,x,y)|w(y)| \Id\mu_g (y) \Id s\in [0,\infty].
$$
Note that $\left\|\cdot\right\|_{g,t}$ becomes an extended norm\footnote{that is, a norm with values in $[0,\infty]$.} on $\IMM(X)$, as the heat kernel is strictly positive. For all $0<t<T$ one has \cite{kuwae2}
\begin{align}\label{kuwe}
\left\|\cdot\right\|_{g,T}\leq \min\{l\in\IN: T< lt\} \left\|\cdot\right\|_{g,t},
\end{align}
which is a consequence of 
\begin{align}\label{mar}
	\int_X p_g(s,x,y) \Id\mu_g(y)\leq 1\quad\text{for all $s>0$, $x\in X$}
\end{align}
and the Chapman-Kolmogorov identity 
$$
p_g(t+s,x,y)=\int_X p_g(t,x,z)p_g(s,z,y)\Id\mu_g(z)\quad\text{for all $t,s>0$, $x,y\in X$.}
$$

\begin{remark} For all $t>0$, $\lambda>0$, $w\in\IMM(X)$ it holds that \cite{batu1}
	\begin{align*}
	&(1-e^{-\lambda t})\sup_{x\in X}\int^{\infty}_0 e^{-\lambda s} \int_X p_g(s,x,y)|w(y)|  \Id\mu_g(y) \Id s\\
	&\leq \left\|w\right\|_{g,t}\\
	&\leq e^{\lambda t}\sup_{x\in X}\int^{\infty}_0 e^{-\lambda s} \int_X p_g(s,x,y)|w(y)|  \Id\mu_g(y) \Id s,
	\end{align*}
	which follows straightforwardly from the Chapman-Kolmogorov identity and \eqref{mar}. In view of
	$$
	(H_g+\lambda)^{-1} w(x)= \int^{\infty}_0 e^{-\lambda s}  \int_X p_g(s,x,y)w(y)  \Id\mu_g(y) \Id s\quad\text{for all $\lambda>0$, $0\leq w\in L^2(X,g)$,}
	$$
	this gives a characterization of the Dynkin norm in terms of the resolvent.
\end{remark}

Let us define extended norms by
$$
\left\| w\right\|_{g,t;A}:=\left\| 1_A w\right\|_{g,t }\quad \text{for every $w\in \IMM(X)$, $A\subset X$ closed, $t >0$.}
$$
Brownian motion\footnote{Note our convention is here that Brownian motion is an $-\Delta$-diffusion rather than an $-\Delta/2$-diffusion.} $((\mathbb{P}^x_{g})_{x\in X}, \mathbb{X})$ on $(X,g)$ is the uniquely determined diffusion on $X$ such that for all $x\in X$, $s>0$ and all Borel sets $B\subset X$ one has 
$$
\mathbb{P}^x_{g}\{\mathbb{X}_{s}\in B, t<\zeta_X\}=\int_B p_g(s,x,y) \Id\mu_g(y),
$$
where $\zeta_U$ denotes the first exit time of $\mathbb{X}$ from an open subset $U$ of the Alexandrov compactification of $X$. It follows that for all $w$, $t$, $A$ as above that 
$$
\left\|w\right\|_{g,t;A}=\sup_{x\in X}\mathbb{E}^x_{g}\left[\int_{0}^t1_{\{s<\zeta_X\}}|1_Aw|(\mathbb{X}_s)\Id s\right].
$$
The right hand side can estimated as follows: with $\sigma_A$ the first hitting time of $\mathbb{\mathbb{X}}$ to $A$ one has 
\begin{align*}
\sup_{x\in X}\mathbb{E}^x_{g}\left[\int_{0}^t1_{\{s<\zeta_X\}}|1_Aw|(\mathbb{X}_s)\Id s\right]&=\sup_{x\in X}\mathbb{E}^x_{g}\left[1_{\{t\geq\sigma_A\}}\int_{\sigma_A}^t1_{\{s<\zeta_X\}}|1_Aw|(\mathbb{X}_s)\Id s\right]\\
&=\sup_{x\in X}\mathbb{E}^x_{g}\left[1_{\{t\geq\sigma_A\}}\int_0^{t-\sigma_A} 1_{\{s<\zeta_X\}} |1_Aw|(\mathbb{X}_{s+\sigma_A})\Id s\right]\\
	&=\sup_{x\in X}\mathbb{E}^x_{g}\left[1_{\{t\geq\sigma_A\}}\mathbb{E}^{\mathbb{X}_{\sigma_A}}_{g}\left[\int_0^{t-\sigma_A}1_{\{s<\zeta_X\}}|1_Aw|(\mathbb{X}_s)\Id s \right]\right] \\
	&\leq \sup_{x\in X}\mathbb{E}^x_{g}\left[1_{\{t\geq\sigma_A\}}\mathbb{E}^{\mathbb{X}_{\sigma_A}}_{g}\left[\int_0^t1_{\{s<\zeta_X\}}|1_Aw|(\mathbb{X}_s)\Id s \right] \right]\\
	&\leq \sup_{x\in A}\mathbb{\mathbb{E}}^x_{g}\left[\int_0^t1_{\{s<\zeta_X\}}|1_Aw|(\mathbb{X}_s)\Id s\right]\\
	&=\sup_{x\in A}\int^t_0\int_A p_{g}(s,x,y) |w(y)| \Id\mu_{g}(y) \Id s, 
\end{align*} 
where we have used the strong Markov property of Brownian motion. In other words, the supremum over $x\in X$ in the definition of $\left\| w\right\|_{g,t;A}$ localizes to $A$. As this observation (which can be seen as some abstract maximum principle) is in the center of this note, we record this fact as:

\begin{theorem}\label{local} For all $t>0$, $w\in\IMM(X)$ and all closed subsets $A\subset X$ one has 
	\begin{align}
		\left\|w\right\|_{g,t;A}= \sup_{x\in A}\int^t_0\int_A p_{g}(s,x,y) |w(y)| \Id\mu_{g}(y) \Id s.
	\end{align}
\end{theorem}

As a first application of this localization we obtain:

\begin{proposition}\label{einb} Let $g$ be a Riemannian metric on $X$. Assume further that $A\subset X$ is closed and that there exists a constant $0<T\leq \infty$ and Borel functions $\phi_1:(0,T)\to [0,\infty)$, $\phi_2:A\to [0,\infty)$ such that for all $0<t<T$, $x,y\in A$ one has
	$$
p_g(t,x,y)\leq \phi_1(t)\phi_2(y).
	$$
	Then for all $w\in \IMM(X)$, $0<t<T$ one has 
	$$
	\left\| w\right\|_{g,t;A}\leq \int^t_0 \phi_1(s)^{\frac{1}{q}}\Id s  \left\| w\right\|_{L^q_{\phi_2}(A,g)}.
	$$
\end{proposition}

\begin{proof} For $q'$ chosen with $1/q'+1/q=1$, $x\in A$, $s>0$ we have 
	\begin{align*}
		&\int_A p_{g}(s,x,y) |w(y)| \Id\mu_{g}(y)\leq \left(\int_A p_{g}(s,x,y)\Id\mu_g(y)\right)^{1/q'}\left(\int_Ap_{g}(s,x,y)|w(y)|^q \Id\mu_{g}(y)\right)^{1/q}\\
		&\leq \phi_1(s)^{\frac{1}{q}}\left(\int_A|w(y)|^q \phi_2(y)\Id\mu_{g}(y)\right)^{1/q},
	\end{align*}
	where we have used \eqref{mar}. Now the claim follows from Theorem \ref{local}.
\end{proof}

\section{Dynkin and Kato class}

In view of \eqref{kuwe} one defines the \emph{Dynkin class} of $(X,g)$ (sometimes also called \emph{extended Kato class} \cite{stollmann,voigt}) by
$$
\widetilde{\mathcal{K}}(X,g):=\left\{w\in\IMM(X): \left\|w\right\|_{g,t}<\infty\quad\text{for some/all $t>0$}\right\}.
$$
The \emph{Kato class} of $(X,g)$ is defined by
$$
\mathcal{K}(X,g):=\left\{w\in\IMM(X): \lim_{t\to 0+}\left\|w\right\|_{g,t}=0\right\}\subset \widetilde{\mathcal{K}}(X,g).
$$
We consider $\widetilde{\mathcal{K}}(X,g)$ as a locally convex space with respect to the family of norms $\left\|\cdot\right\|_{g,t}$, $t>0$. Clearly, $\widetilde{\mathcal{K}}(X,g)$ is countably normed (by letting $t$ run through $\IN$).

\begin{remark} Again by \eqref{kuwe}, the topology on $\widetilde{\mathcal{K}}(X,g)$ is in fact equal to the topology induced by the single norm $\left\|w\right\|_{g,t=1}$. However, we believe it is more natural to consider all $\left\|\cdot\right\|_{g,t}$.
\end{remark}

Using \eqref{mar} one finds the trivial continuous inclusion 
\begin{align}\label{inc}
L^\infty(X,g)\subset \mathcal{K}(X,g),\quad\text{with $\left\|w\right\|_{g,t}\leq t \left\|w\right\|_{L^\infty(X,g)}$ for all $t>0$, $w\in L^\infty(X,g)$.} 
\end{align}

\begin{proposition}\label{global} Let $g$ be a Riemannian metric on $X$. Then $\widetilde{\mathcal{K}}(X,g)$ is a Fr\'{e}chet space and $\mathcal{K}(X,g)$ is a closed subspace of $\widetilde{\mathcal{K}}(X,g)$ (thus also a Fr\'{e}chet space).
\end{proposition}

\begin{proof} To see that $\widetilde{\mathcal{K}}(X,g)$ is complete, let us first record that for all $w\in\IMM(X)$, $t>0$ and all compact sets $K\subset X$ one has the simple estimate
	\begin{align}\label{emb}
		\int_K  |w(y)|  \Id\mu_g(y)\leq \frac{2}{t}\left(\min_{s\in [t/2,t], z\in K ,y\in K}p_g(s,z,y)\right)^{-1}\sup_{x\in K}\int^t_0\int_K p_g(s,x,y) |w(y)| \Id s \Id\mu_g(y),
	\end{align}
	where 
	$$
	  \min_{s\in [t/2,t], z\in K ,y\in K}p_g(s,z,y)>0
	$$
	follows from $p_g>0$. Now let $w_n$ be a Cauchy sequence in $\widetilde{\mathcal{K}}(X,g)$. Then by \eqref{emb} the sequence $w_n$ has a subsequence $w_{n_j}$ which converges almost everywhere to some equivalence class of Borel functions $w$ on $X$. Given $t,\epsilon>0$ pick $N=N(\epsilon,t)\in\IN$ such for all $j,l\geq N$ one has $\left\|w_{n_j}-w_l\right\|_{g,t}<\epsilon/2$. For such $l$ one has 
	\begin{align*}
		&\left\|w-w_l\right\|_{g,t}= 	\sup_{x\in X}   \int^t_0\int_X p_g(s,x,y) |w-w_l|(y) \Id s \Id\mu_g(y)\\
		& \leq  \liminf_j \sup_{x\in X}  \int^t_0\int_X p_g(s,x,y) |w_{n_j}-w_l|(y) \Id s  \Id\mu_g(y)= \liminf_j \left\|w_{n_j}-w_l\right\|_{g,t} < \epsilon,	
	\end{align*}	
	where we have used Fatou's lemma. This shows that $w\in \widetilde{\mathcal{K}}(X,g)$ and $w_l\to w$ in $\widetilde{\mathcal{K}}(X,g)$.\\
	To see that $\mathcal{K}(X,g)$ is a closed subspace, let $w_n$ be a sequence in $\mathcal{K}(X,g)$ which converges to $w\in \widetilde{\mathcal{K}}(X,g)$. Let $t,\epsilon>0$, pick $n_0=n_0(\epsilon)\in\IN$ with $\left\|w-w_{n_0}\right\|_{g,t=1}<\epsilon/2$ and pick $t_0=t_0(\epsilon,n_0)\in (0,1]$ such that for all $0<t<t_0$ one has $\left\|w_{n_0}\right\|_{g,t}<\epsilon/2$. Then for these $t$'s one has 
	$$
	\left\|w\right\|_{g,t}\leq \left\|w-w_{n_0}\right\|_{g,t=1}+\left\|w_{n_0}\right\|_{g,t}<\epsilon.
	$$
	This completes the proof.
\end{proof}

Assume we are given a metric Dirichlet space $(\mathsf{X},\mathsf{d}, \mathsf{m},\mathsf{E})$. With $\mathsf{H}\geq 0$ the self-adjoint operator in $L^2(\mathsf{X},\mathsf{m})$ associated with the Dirichlet form $\mathsf{E}$, we assume that the heat semigroup of $\mathsf{H}$ has a jointly continuous integral kernel
$$
\mathsf{p}:(0,\infty)\times \mathsf{X}\times \mathsf{X}\longrightarrow [0,\infty).
$$
Let $(x,r)\mapsto \mathsf{m}(x,r)$ denote the volume function on $(\mathsf{X},\mathsf{d}, \mathsf{m})$, that is, $\mathsf{m}(x,r)$ is the volume of an open ball with radius $r$ around $x$. Then $(\mathsf{X},\mathsf{d}, \mathsf{m},\mathsf{E})$ is said to satisfy 

\begin{itemize}
\item the \emph{local-in-time Gaussian heat kernel upper estimate}, if there exist constants $\alpha, \beta,\gamma>0$ such that for all $t>0$, $x,y\in \mathsf{X}$, 
\begin{align*}
\bf{GUE_\loc}:\quad\quad \quad \mathsf{p}(t,x,y)\leq 		\alpha \mathsf{m}(x,\sqrt{t})^{-1}e^{-\beta\frac{\mathsf{d}(x,y)^2}{t}} e^{\gamma t},
\end{align*}
\item the \emph{local-in-time Gaussian heat kernel lower estimate}, if there exist constants $\alpha', \beta',\gamma'>0$ such that for all $t>0$, $x,y\in  \mathsf{X}$, 
\begin{align*}
\bf{GLE_\loc}:\quad \quad\quad \mathsf{p}(t,x,y)\geq 		\alpha' \mathsf{m}(x,\sqrt{t})^{-1}e^{-\beta'\frac{\mathsf{d}(x,y)^2}{s}} e^{-\gamma' t},
\end{align*}
\item \emph{locally uniformly $N$-volume doubling} for a given $N>0$, if there exists a constant $a>0$ such that for all $0<r<R$, $x\in \mathsf{X}$,
\begin{align*}
\bf{VD_\loc}(N):\quad \quad \mathsf{m}(x,R)\leq a e^{aR}\left(\frac{R}{r}\right)^N\mathsf{m}(x,r). 
\end{align*}
\end{itemize}

If $(X,\varrho_g,\mu_g)$ satisfies $\bf{GUE_\loc}$ and $\bf{VD_\loc}(N)$ for some $N\geq 2$, then one has a continuous inclusion
\begin{align}\label{ind}
L^q_{\phi_g}(X,g)+L^\infty(X,g)\subset \mathcal{K}(X,g)\quad\text{for all $q\in (N/2,\infty)$.}
\end{align}
Indeed, $\bf{VD_\loc}(N)$ implies that for all $0<t<1$ one has
$$
\mu_g(x,\sqrt{t})^{-1}\leq a e^{a}t^{-N_g/2}\mu_g(x,1)^{-1},
$$	
so that $\bf{GUE_\loc}$ implies
$$
p_g(t,x,y)\leq \alpha a e^{a} t^{-N_g/2}\mu_g(x,1)^{-1}e^{-\beta\frac{\varrho_g(x,y)^2}{t}} \leq \alpha a e^{a} t^{-N_g/2}\mu_g(x,1)^{-1},
$$
and (\ref{ind}) follows from Proposition \ref{einb} in combination with \eqref{inc}.\vspace{1mm}

As a consequence of this observation and the recently established results from \cite{carron3} we would like to highlight:

\begin{theorem}\label{main2} Let $m\geq 2$, let $g$ be a complete Riemannian metric on $X$ and define
	\begin{align*}
	&\sigma_g:X\longrightarrow \IR,\quad \sigma_g(x):= \text{smallest eigenvalue of $\mathrm{Ric}_g(x)\in \mathrm{End}(T_x X)$},\\
	&\sigma_g^-:X\longrightarrow [0,\infty),\quad \sigma_g^-(x):= -\min(\sigma^g(x) , 0),\\
    &\phi_g:X\longrightarrow (0,\infty),\quad\phi_g(x):=1/\mu_g(x,1),\\
    &t^*_g:=\inf\left\{t>0: \left\|w\right\|_{g,t}<1/(3m-6)\right\}\in [0,\infty],\\
    &N_g:=m+4(m-2)^2t^*_g     \in [2,\infty].
\end{align*}
\emph{a)} If $\sigma_g^-	\in  \mathcal{K}(X,g)$, then $(X,\varrho_g,\mu_g)$ satisfies $\bf{GUE_\loc}$, $\bf{GLE_\loc}$ and $\bf{VD_\loc}(N_g)$, in particular, for all $q\in (N_g/2,\infty)$ there is a continuous inclusion $L^q_{\phi_g}(X,g)+L^\infty(X,g)\subset \mathcal{K}(X,g)$.\vspace{1mm}

\emph{b)} Assume that $\sigma_g^-\in L^q_{\phi_g}(X,g)+L^\infty(X,g)$ for some $q\in (1,\infty)$. Then the following statements are equivalent:
\begin{itemize}
	\item $\sigma_g^-\in \mathcal{K}(X,g)$ and $q>N_g/2$.
	\item $(X,\varrho_g,\mu_g)$ satisfies $\bf{GUE_\loc}$ and $\bf{VD_\loc}(N)$ for some $N\geq 2$ with $q>N/2$ .
\end{itemize}

\end{theorem}

\begin{proof} a) By Theorem A in \cite{carron3} the metric measure space $(X,\varrho_g,\mu_g)$ is quasi-isometric to a (smooth) $\mathsf{RCD}(K_g,N_g)$ space for some $K_g\in\IR$. On the other hand, $\mathsf{RCD}(K_g,N_g)$ spaces (with their Cheeger form) satisfy $\bf{GUE_\loc}$, $\bf{GLE_\loc}$ and $\bf{VD_\loc}(N_g)$ (cf. \cite{jiang}). Moreover, $\bf{GUE_\loc}$ together with $\bf{VD_\loc}(N_g)$, as well as $\bf{GUE_\loc}$ together with $\bf{GLE_\loc}$ and $\bf{VD_\loc}(N_g)$ are stable under quasi-isometry. It follows that $(X,\varrho_g,\mu_g)$ satisfies $\bf{GUE_\loc}$, $\bf{GLE_\loc}$ and $\bf{VD_\loc}(N_g)$. \vspace{1mm}
	
b) This is an immediate consequence of part a) and (\ref{ind}).

\end{proof}

\section{Local Kato class and local Dynkin class}

We start this section with the following application of Theorem \ref{local}:

\begin{theorem}\label{main3} \emph{a)} Let $g_1,g_2$ be Riemannian metrics on $X$. Then for all compact $K\subset X$ there exist constants $T_K=T_K(g_1,g_2)>0$ and $C_K=C_K(g_1,g_2)>1$ such that for all $0<t<T_K$, $w\in \IMM(X)$, 
$$
(1/C_K)\left\| w\right\|_{g_2,t;K}\leq \left\| w\right\|_{g_1,t;K}\leq C_K\left\| w\right\|_{g_2,t;K}.
$$
\emph{b)} Assume $m\geq 2$, $q\in (m/2,\infty)$ and let $g$ be a Riemannian metric on $X$. Then for every $K\subset X$ compact there exist constants $T_K=T_K(g)>0$ and $C_{K,q,m}=C_{K,q,m}(g)>0$ such that for all $0<t<T_K$ , $w\in \IMM(X)$,
$$
\left\| w\right\|_{g,t;K}\leq C_{K,q,m}\int^t_0 s^{-\frac{m}{2q}}\Id s \left\|w\right\|_{L^q(K,g)}.
$$
\end{theorem}

\begin{proof} Let us first show for every Riemannian metric $g$ on $X$ and every $x_0\in X$, there exists an open relatively compact neighbourhood $V_{x_0}=V_{x_0}(g)$ of $x_0$ and numbers 
\begin{align*}
&S_{x_0}=S_{x_0}(g)>0, \quad \gamma_{x_0}=\gamma_{x_0}(g)>0,\quad  \gamma'_{x_0}=\gamma'_{x_0}(g)>0,\\
&\delta_{x_0}=\delta_{x_0}(g)>0,\quad \delta'_{x_0}=\delta'_{x_0}(g)>0,
\end{align*}
such that for all $0<s<S_{x_0}$, $x,y\in \overline{V_{x_0}}$, 
	\begin{align}\label{edc}
		\gamma'_{x_0}s^{-m/2}e^{-\gamma_{x_0}\frac{\varrho_g(x,y)^2}{s}}\leq 	p_g(s,x,y)\leq 	\delta'_{x_0}s^{-m/2}e^{-\delta_{x_0}\frac{\varrho_g(x,y)^2}{s}}.
	\end{align}
To see this, we record that by \cite{sturm} there exist constants 
$$
\gamma_{x_0}=\gamma_{x_0}(g)>0,\quad  \gamma''_{x_0}=\gamma''_{x_0}(g)>0,\quad \delta_{x_0}=\delta_{x_0}(g)>0,\quad \delta''_{x_0}=\delta''_{x_0}(g)>0
$$	
such that for all $0<s<S_{x_0}:=\min(\mathrm{inj}_g(x_0),1)^2$ and all $x,y\in \overline{B}_g(x_0,r)$,
	\begin{align*}
		\gamma_{x_0,r}''\mu_g(x,\sqrt{s})^{-1}e^{-\gamma_{x_0,r}\frac{\varrho_g(x,y)^2}{s}}\leq p_g(s,x,y)\leq 	\delta_{x_0,r}''\mu_g(x,\sqrt{s})^{-1}e^{-\delta_{x_0,r}\frac{\varrho_g(x,y)^2}{s}}.
\end{align*}

Using that there exist $\delta'''_{x_0,r}=\delta'''_{x_0,r}(g)>0$, $\gamma'''_{x_0,r}=\gamma'''_{x_0,r}(g)>0$ such that for all $0<s<r^2$,
	$$
	\delta'''_{x_0,r} s^{-m/2}\leq \mu_g(x,\sqrt{s})^{-1}\leq \gamma'''_{x_0,r} s^{-m/2},
	$$
	a consequence of local volume doubling and local Bishop-Gromov volume estimate, we get \eqref{edc}.\\
Now with the open relatively compact set $U_{x_0}:=V_{x_0}(g_1)\cap V_{x_0}(g_2)$ one gets $a_{x_0}=a_{x_0}(g_1,g_2)>1$ such that $(1/a_{x_0})g_2\leq g_1\leq a_{x_0}g_2$ on $\overline{U_{x_0}}$, thus we also get constants
\begin{align*}
&b_{x_0}=b_{x_0}(g_1,g_2)>1,\quad c_{x_0}=c_{x_0}(g_1,g_2)>1\quad\text{such that}\\
	&(1/b_{x_0})\varrho_{g_2}\leq \varrho_{g_1} \leq b_{x_0}\varrho_{g_2}\quad\text{on $\overline{U_{x_0}}\times \overline{U_{x_0}}$},\\
& (1/c_{x_0})\mu_{g_2}\leq \mu_{g_1} \leq c_{x_0}\mu_{g_2}, \quad\text{on Borel subsets of $\overline{U_{x_0}}$}.
\end{align*}
It follows then from \eqref{edc} that there exists $C_{x_0}=C_{x_0}(g_1,g_2)>1$ such that with $T_{x_0}:=\min(S_{x_0}(g_1), S_{x_0}(g_2))$ for all $0<t<T_{x_0}$ one has
\begin{align*}
	&(1/C_{x_0})\sup_{x\in \overline{U_{x_0}}}\int^t_0\int_{\overline{U_{x_0}}} p_{g_2}(s,x,y) |w(y)| \Id\mu_{g_2}(y) \Id s  \\
		&\leq \sup_{x\in \overline{U_{x_0}}}\int^t_0\int_{\overline{U_{x_0}}} p_{g_1}(s,x,y) |w(y)| \Id\mu_{g_1}(y) \Id s  \\
		&\leq C_{x_0}\sup_{x\in \overline{U_{x_0}}}\int^t_0\int_{\overline{U_{x_0}}} p_{g_2}(s,x,y) |w(y)| \Id\mu_{g_2}(y) \Id s.
	\end{align*}	
Theorem \ref{local} then shows
$$
(1/C_{x_0})\left\| w\right\|_{g_2,t;\overline{U_{x_0}}}\leq \left\| w\right\|_{g_1,t;\overline{U_{x_0}}}\leq C_{x_0}\left\| w\right\|_{g_2,t;\overline{U_{x_0}}  }.
$$
The claim now follows from a simple covering argument.\vspace{1mm}

b) In view of \eqref{edc}, given $x_0 \in K$ we can apply Proposition \ref{einb} to get that for all $0<t<T_{x_0}$ , $w\in \IMM(X)$,
$$
\left\| w\right\|_{g,t;\overline{U_{x_0}}}\leq C_{x_0,q,m}\int^t_0 s^{-\frac{m}{2q}}\Id s \left\|w\right\|_{L^q(\overline{U_{x_0}},g)},
$$
where $C_{x_0,q,m}=C_{x_0,q,m}(g)>0$, $T_{x_0}=T_{x_0}(g)>0$. Again a simple covering argument yields the claim.
\end{proof}

Define 
$$
\widetilde{\mathcal{K}}_\loc(X,g):=\{w\in\IMM(X): 1_K w\in \widetilde{\mathcal{K}}(X,g)\>\text{for all compact $K\subset X$}\},
$$
and 
$$
\mathcal{K}_\loc(X,g):=\{w\in\IMM(X): 1_K w\in \mathcal{K}(X,g)\>\text{for all compact $K\subset X$}\}\subset \widetilde{\mathcal{K}}_\loc(X,g).
$$

Then $\widetilde{\mathcal{K}}_\loc(X,g)$ becomes a locally convex space with respect to the family of norms
$$
\left\| \cdot\right\|_{g,t;K},\quad \text{$t>0$, $K\subset X$ compact,}
$$
which is in fact countably normed (take a countable exhaustion of $X$ with compacts and let $t$ run through $\IN$) and thus is metrizable. It is checked precisely as in the proof of Proposition \ref{global} that $\widetilde{\mathcal{K}}_\loc(X,g)$ is a Fr\'{e}chet space and that $\mathcal{K}_\loc(X,g)$ is a closed subspace of $\widetilde{\mathcal{K}}_\loc(X,g)$.\vspace{1mm}

From the results obtained in this section, we immediately get that the local Dynkin resp. the local Kato class is defined on any manifold:

\begin{propandef} \emph{a)} For any two Riemannian metrics $g_1,g_2$ on $X$ we have $\widetilde{\mathcal{K}}_\loc(X,g_1)=\widetilde{\mathcal{K}}_\loc(X,g_2)$ as locally convex spaces, so that we get a well-defined Fr\'{e}chet space $\widetilde{\mathcal{K}}_\loc(X)$ called the \emph{local Dynkin class on $X$}. Likewise, we have $\mathcal{K}_\loc(X,g_1)=\mathcal{K}_\loc(X,g_2)$ as locally convex spaces, so that we get the closed subspace $\mathcal{K}_\loc(X)\subset \widetilde{\mathcal{K}}_\loc(X)$ called the \emph{local Kato class on $X$}.\\
\emph{b)} One has continuous inclusions $L^\infty_\loc(X)\subset \mathcal{K}_\loc(X)$ and $\widetilde{\mathcal{K}}_\loc(X)\subset L^1_\loc(X)$, and if $m\geq 2$, $q\in (m/2,\infty)$ then we also have a continuous inclusion $L^q_\loc(X)\subset \mathcal{K}_\loc(X)$.
\end{propandef}

\begin{proof} a) This follows immediately from Theorem \ref{main3} a) and \eqref{kuwe}.\\
b) $L^\infty_\loc(X)\subset \mathcal{K}_\loc(X)$ is a consequence of \eqref{inc}, $\widetilde{\mathcal{K}}_\loc(X)\subset L^1_\loc(X)$ is a consequence of estimate \eqref{emb}, and $L^q_\loc(X)\subset \mathcal{K}_\loc(X)$ follows from Theorem \ref{main3} b). 
\end{proof}

We remark that the inclusion $L^q_\loc(X)\subset \mathcal{K}_\loc(X)$ has been shown \cite{batu3} using the parabolic $L^1$-mean value inequality.


\begin{proposition}\label{dens} $C^\infty_c(X)$ is dense in $\mathcal{K}_\loc(X)$.
\end{proposition}

\begin{proof} Pick a Riemannian metric $g$ on $X$ so that $\mathcal{K}_\loc(X)=\mathcal{K}_\loc(X,g)$. Let us first show that compactly supported elements of $\mathcal{K}_\loc(X,g)$ are dense in $\mathcal{K}_\loc(X,g)$. Indeed, given $w\in \mathcal{K}_\loc(X,g)$, let $K_n$ be an exhaustion of $X$ with compact subsets and set $w_n:=1_{K_n}w\in \mathcal{K}_\loc(X,g)$. Then, given $K\subset X$ compact, picking $n_0=n_0(K)\in\IN$ with $K\subset K_n$ for all $n\geq n_0$ we have $1_K1_{K_n}=1_K$, and for those $n$ it follows that
$$
\left\| w-w_n\right\|_{g,t;K}=\sup_{x\in X}\int^t_0\int_{X} p_g(s,x,y) |1_Kw(y)-1_K1_{K_n}w(y)| \Id y \Id s=0.
$$	
To proceed further, for all compact subsets $K\subset X$ which are included in a chart of $X$ and all $t>0$, $w\in\IMM(X)$ we set
	$$
	\left\| w\right\|_{t;K}:=\sup_{x\in K}\int^t_0s^{-m/2}\int_{K} e^{-\frac{|x-y|^2}{4s}} |w(y)| \Id y \Id s.
	$$
By the proof of Theorem \ref{main3} a), for every point $x_0\in X$ there exists a chart $V_{x_0}=V_{x_0}(g)$ around $x_0$, and numbers $C_{x_0}=C_{x_0}(g)>1$ and $T_{x_0}=T_{x_0}(g)>0$ such that for all $0<t<T_{x_0}$, all compact $K\subset V_{x_0}$ and all $w\in \IMM(X)$ one has
\begin{align}\label{ayaa}
(1/C_{x_0})	\left\| w\right\|_{t;K}\leq 	\left\| w\right\|_{g,t;K}\leq 	C_{x_0}\left\| w\right\|_{t;K}.
\end{align}
Assume $w\in \mathcal{K}_\loc(X,g)$ is compactly supported pick finitely many $x_1,\dots, x_l$ in $\mathrm{supp}(w)$ such that $\mathrm{supp}(w)\subset \bigcup_{l=1}V_{x_j}$. Pick a partition of unity $\psi_j\in C^{\infty}_c(V_{x_j})$ subordinate to $V_{x_j}$. Then $w_j:=\psi_jw$ is in $\mathcal{K}_\loc(X,g)$ and by \eqref{ayaa} and Proposition 2.3 in \cite{hundert}, which states that $C^\infty_{c}(\IR^m)$ is dense in $\mathcal{K}_\loc(\IR^m,g_{\mathrm{Eucl}})$, we can pick a sequence of $w_{j,n}\in C^\infty_{c}(X)$, $n\in\IN$, such that $w_{j,n}\to w_j$ as $n\to\infty$ in $\mathcal{K}_\loc(X,g)$. It follows that $w_n:=\sum^l_{j=1}w_{j,n}\in C^\infty_{c}(X)$ converges to $w$ in $\mathcal{K}_\loc(X,g)$.  
\end{proof}

\section{Continuity properties of Schr\"odinger semigroups}

Let $g$ be a Riemannian metric on $X$ and let $w\in\IMM(X)$ be a \emph{Kato decomposable} potential on $(X,g)$, which means that $w=w^+-w^-$ for some $0\leq w^\pm \in\IMM(X)$ with $w^+\in \mathcal{K}_\loc(X)$ and $w^-\in\mathcal{K}(X,g)$. It follows that $w\in L^1_\loc(X)$ and moreover that the symmetric bilinear form in $L^2(X,g)$ defined by 
$$
C_c^\infty(X)\times C_c^\infty(X)\ni (\psi_1,\psi_2)\longmapsto \int_X (-\Delta_g\psi_1)\cdot \psi_2 \Id\mu_g+  \int_X w \cdot \psi_1 \cdot \psi_2  \Id\mu_g\in \IR
$$ 
is closable and semibounded (from below) \cite{batu1, batu4}. We denote the semibounded self-adjoint operator in $L^2(X,g)$ induced by the closure of this form by $H^w_g$. We are going to prove the following result, which generalizes Theorem X.1 in \cite{batu1} in the sense that not global assumption on $w^+$ is required. We refer the reader also to \cite{hundert} for the state of the art in Euclidean space concerning continuity properties of Schr\"odinger semigroups.

\begin{theorem}\label{cou} In the above situation, for all $\psi\in L^2(X,g)$ the assignment $(t,x)\mapsto e^{-t H^w_g}\psi(x)$ has a jointly continuous representative on $(0,\infty)\times X$.
\end{theorem}

We record the following well-known fact, usually referred to as \emph{Khashminski's lemma}, in a form which is suitable for the proof:

\begin{lemma}Let $0\leq w\in \IMM(X)$, let $U\subset X$ be open and set $w^U:=w|_{U}$, $g_U:=g|_{U}$. Assume there exists $t_0>0$ with $\left\| w^U\right\|_{g_U,t_0}<1$. Then with 
	$$
	C(g,U,w,t_0):=\frac{1}{1-\left\| w^U\right\|_{g_U,t_0}},
	$$
	for all $t>0$, $x\in U$ one has 
	$$
	\mathbb{E}^x_g\left[1_{\{t<\zeta_U\}}e^{\int^t_0w (\mathbb{X}_s) \Id s}\right]\leq C(g,U,w,t_0)\exp\left(\frac{t}{t_0}\log\big(C(g,U,w,t_0)\big)\right)<\infty.
	$$
\end{lemma}

A proof of this result can be found in \cite{batu1} (just apply Lemma VI.8 therein to the Riemannian manifold $(U,g|_U)$).

\begin{proof}[Proof of Theorem \ref{cou}] Let $U_l$, $l\in\IN$, be an exhaustion of $X$ with open relatively compact subsets. We set 
$$
g_l:=g|_{U_l},\quad \zeta_l:=\zeta_{U_l},\quad \zeta:=\zeta_X,\quad \psi^{(l)}:=\psi|_{U_l},\quad w^{(l)}:=w|_{U_l}.
$$
Because of 
$$
p_{g_l}(x,y,t)\leq p_{g}(x,y,t) \quad\text{for all $t>0$, $x,y\in U_l$,}
$$
we have $w^{(l)}\in\mathcal{K}(U_l,g_l)$ as $U_l$ is relatively compact, and so $H^{w}_{g,l}:=H^{w_l}_{g_{l}}$ is well-defined in $L^2(U_l,g_l)$. In view of the Feynman-Kac formulae for $e^{-t H^{w}_{g_l}}$ and $e^{-t H^{w}_{g} }$ (which hold in this generality , c.f. \cite{bei}), a pointwise defined representative of $(t,x)\mapsto e^{-t H^{w}_{g_l}   }\psi^{(l)}(x)$ on $(0,\infty)\times U_l$ is defined by
	$$
	e^{-t H^{w}_{g,l}   }\psi^{(l)}(x)=\mathbb{E}^x_g\left[1_{\{t<\zeta_l\}}e^{-\int^t_0w (\mathbb{X}_s) \Id s}\psi(\mathbb{X}_t)\right],
	$$
and a pointwise defined representative of $(t,x)\mapsto e^{-t H^{w}_{g}   }\psi(x)$ on $(0,\infty)\times X$ is defined by
$$
e^{-t H^{w}_{g}   }\psi(x)=\mathbb{E}^x_g\left[1_{\{t<\zeta\}}e^{-\int^t_0w (\mathbb{X}_s) \Id s}\psi(\mathbb{X}_t)\right].
$$	
We pick a sequence $w_n$ in $C^\infty_c(X)$ with $w_n\to w$ in $\mathcal{K}_\loc(X)$.

Step 1: For all $\alpha\geq 0$ and all compact intervalls $I\subset  (0,\infty)$ one has
$$
\sup_{x\in U_l,n\in\IN, t\in I}\mathbb{E}^x_g\left[1_{\{t<\zeta_l\}}e^{\alpha\int^t_0|w_n (\mathbb{X}_s)| \Id s}\right]<\infty.
$$
Proof: Upon replacing $w$ with $\alpha w$ we can assume $\alpha=1$. Clearly, we have $w^{(l)}_n\to w^{(l)}$ in $\mathcal{K}(U_l,g_l)$ as $n\to\infty$ as $U_l$ is relatively compact. It follows that 
$$
\left\| w^{(l)}_n\right\|_{g_l,t}\leq \left\| w^{(l)}-w^{(l)}_n\right\|_{g_l,t=1}+\left\| w^{(l)}\right\|_{g_l,t},
$$
and so there exists $n_0\in\IN$ and $t_0>0$ such that for all $n\geq n_0$,
$$
\left\| w^{(l)}_n\right\|_{g_l,t_0}< 1,
$$ 
and the claim follows from Khashminski's lemma. \vspace{1mm}

Step 2: For all $\psi\in C^\infty_c(X)$ the function $(t,x)\mapsto e^{-t H^{w}_{g,l}   }\psi^{(l)}(x)$ is jointly continuous on $(0,\infty)\times U_l$. 
Proof: The map $(t,x)\mapsto e^{-t H^{w_n}_{g,l}   }\psi^{(l)}(x)$ is smooth by local parabolic regularity. For all compact intervals $I\subset (0,\infty)$, $t\in I$, $x\in U_l$, using the Feynman-Kac formulae and w.l.o.g. $\left\|\psi\right\|_\infty\leq 1$ we get
$$
	\left|e^{-t H^{w}_{g,l}   }\psi^{(l)}(x)-e^{-t H^{w_n}_{g,l}   }\psi^{(l)}(x)\right|\leq  \mathbb{E}^x_g\left[1_{\{t<\zeta_l\}}\left|e^{-\int^t_0w (\mathbb{X}_s) \Id s}-e^{-\int^t_0w_n (\mathbb{X}_s) \Id s}\right|\right].
$$	
Using the inequality 
$$
|e^{-a}-e^{-b}|\leq  e^{1+a^-+b^-}|a-b|^{1/2},\quad\text{for all $a,b\in\IR$},
$$
the latter is
\begin{align*}
&\leq e \mathbb{E}^x_g\left[1_{\{t<\zeta_l\}}e^{\int^t_0w_n^- (\mathbb{X}_s) \Id s}e^{\int^t_0w^- (\mathbb{X}_s) \Id s}\left|\int^t_0w(\mathbb{X}_s)-w_n (\mathbb{X}_s) \Id s\right|^{\frac{1}{2}}\right]\\
&\leq e \mathbb{E}^x_g\left[1_{\{t<\zeta_l\}}e^{2\int^t_0w_n^- (\mathbb{X}_s) \Id s}e^{2\int^t_0w^- (\mathbb{X}_s) \Id s}\right]^{\frac{1}{2}}\mathbb{E}^x_g\left[1_{\{t<\zeta_l\}}\left|\int^t_0w(\mathbb{X}_s)-w_n (\mathbb{X}_s) \Id s\right|\right]^{\frac{1}{2}}\\
&\leq e \mathbb{E}^x_g\left[1_{\{t<\zeta_l\}}e^{4\int^t_0w_n^- (\mathbb{X}_s) \Id s}\right]^{\frac{1}{4}}\mathbb{E}^x_g\left[1_{\{t<\zeta_l\}}e^{4\int^t_0w^- (\mathbb{X}_s) \Id s}\right]^{\frac{1}{4}}\mathbb{E}^x_g\left[1_{\{t<\zeta_l\}}\int^t_0\left|w(\mathbb{X}_s)-w_n (\mathbb{X}_s) \right|\Id s\right]^{\frac{1}{2}}\\
&\leq C_l \mathbb{E}^x_g\left[1_{\{t<\zeta_l\}}\int^t_0\left|w(\mathbb{X}_s)-w_n (\mathbb{X}_s) \right|\Id s\right]^{\frac{1}{2}},
    \end{align*}
where
$$
C_l:=e \sup_{x\in U_l,n\in\IN, t\in I} \left(\mathbb{E}^x_g\left[1_{\{t<\zeta_l\}}e^{4\int^t_0w_n^- (\mathbb{X}_s) \Id s}\right]^{\frac{1}{4}}\mathbb{E}^x_g\left[1_{\{t<\zeta_l\}}e^{4\int^t_0w^- (\mathbb{X}_s) \Id s}\right]^{\frac{1}{4}}\right)<\infty
$$
follows from Khashminski's lemma and Step 1. Finally, we have
\begin{align*}
&\mathbb{E}^x_g\left[1_{\{t<\zeta_l\}}\int^t_0\left|w(\mathbb{X}_s)-w_n (\mathbb{X}_s) \right|\Id s\right],\\
&\leq  \mathbb{E}^x_g\left[\int^t_01_{\{s<\zeta_l\}}\left|w(\mathbb{X}_s)-w_n (\mathbb{X}_s) \right|\Id s\right]\\
 &\leq    \int^{\max I}_0 \int_Xp_{g}(x,y,s)\left|1_{\overline{U_l}}(y)w(y)-1_{\overline{U_l}}(y)w_n (y)\right|\Id\mu_g(y)\Id s\\
&\leq  \left\|w_n-w\right\|_{g,\max I;\overline{U_l}},
\end{align*}
which proves Step 2.\vspace{1mm}

Step 3: For all $\psi\in C^\infty_c(X)$ the function $(t,x)\mapsto e^{-t H^{w}_{g}   }\psi(x)$ is jointly continuous on $(0,\infty)\times X$. \\
	Proof: Fix compacts $I\subset (0,\infty)$, $K\subset X$ and let $l$ be large enough with $K\subset U_l$. For all $x\in K$, $t\in I$ we have 
	\begin{align*}
	&\left|e^{-t H^{w}_{g}   }\psi(x)-e^{-t H^{w}_{g,l}   }\psi^{(l)}(x)\right|\leq \left\|\psi\right\|_\infty\mathbb{E}^x_g\left[\left(1_{\{t<\zeta\}}-1_{\{t<\zeta_l\}}\right)e^{-\int^t_0w (\mathbb{X}_s)\Id s }\right]\\
	&\leq \left\|\psi\right\|_\infty\mathbb{E}^x_g\left[1_{\{t<\zeta\}}-1_{\{t<\zeta_l\}}\right]^{1/2}\mathbb{E}^x_g\left[1_{\{t<\zeta\}}e^{-2\int^t_0w (\mathbb{X}_s)\Id s }\right]^{1/2}\\
    &\leq \mathbb{E}^x_g\left[1_{\{t<\zeta\}}-1_{\{t<\zeta_l\}}\right]^{1/2} \left\|\psi\right\|_\infty \sup_{t\in I, x\in X}\mathbb{E}^x_g\left[1_{\{t<\zeta\}}e^{-2\int^t_0w (\mathbb{X}_s)\Id s }\right]^{1/2},
\end{align*}
where 
$$
\sup_{t\in I, x\in X}\mathbb{E}^x_g\left[1_{\{t<\zeta\}}e^{-2\int^t_0w (\mathbb{X}_s) \Id s}\right]<\infty
$$
follows from Khashminski's lemma. Finally, 
$$
\mathbb{E}^x_g\left[1_{\{t<\zeta\}}-1_{\{t<\zeta_l\}}\right]= \int_X p_g(t,x,y) d\mu_g(y)- \int_{U_l} p_{g_l}(t,x,y) d\mu_g(y),
$$	 
which converges uniformly in $t\in I$, $x\in K$ (in fact, in $C^\infty$ in any open nbh of $K\times I$) to $0$ as $l\to\infty$. Thus the proof of Step 3 is completed by Step 2.\vspace{1mm} 
	
Step 4: For all $\psi\in L^2(X,g)$ the function $(t,x)\mapsto e^{-t H^{w}_{g}   }\psi(x)$ is jointly continuous on $(0,\infty)\times X$.\\
Proof: Pick a sequence $\psi_n$ in $ C^\infty_c(X)$ with $\psi_n\to\psi$ in $L^2(X,g)$ as $n\to\infty$. Let $K\subset X$ and $I\subset (0,\infty)$ both be compact. Then for all $x\in X$, $t\in I$,
\begin{align*}
&\left|e^{-t H^{w}_{g}   }\psi(x)-e^{-t H^{w}_{g}   }\psi_n(x)\right|\leq \mathbb{E}^x_g\left[1_{\{t<\zeta\}}e^{-\int^t_0w (\mathbb{X}_s) \Id s }\right]^{1/2} \mathbb{E}^x_g\left[1_{\{t<\zeta\}}|\psi_n(\mathbb{X}_t)-\psi(\mathbb{X}_t)|^2\right]^{1/2}\\
&\leq\left(\int_X p_g(t,x,y)|\psi_n(y)-\psi(y)|^2 \Id\mu_g(y)\right)^{1/2}\sup_{t\in I, x\in K}\mathbb{E}^x_g\left[1_{\{t<\zeta\}}e^{-\int^t_0w (\mathbb{X}_s) \Id s }\right]^{1/2}\\
&\leq \left(\int_X |\psi_n(y)-\psi(y)|^2 \Id\mu_g(y)\right)^{1/2}\left(\sup_{x\in K,y\in X,t\in I} p_g(t,x,y)\right)^{1/2}\sup_{t\in I, x\in K}\mathbb{E}^x_g\left[1_{\{t<\zeta\}}e^{-\int^t_0w (\mathbb{X}_s) \Id s }\right]^{1/2},
\end{align*}	
where
$$
\sup_{x\in K,y\in X,t\in I} p_g(t,x,y)<\infty
$$
follows from Theorem 2.9 in \cite{batu3} in combination with Example 2.6.2 therein. So the assertion follows from Step 3, and this completes the proof. 
\end{proof}

%




\emph{Acknowledgements.} The authors would like Gilles Carron for pointing out a mistake in a previous version of Theorem 3.3, and J\"urgen Voigt for very helpful discussions. The first named author has been financially supported by the DFG project SPP 2026: Geometry at Infinity. 
The second named author has been supported in part by JSPS Grant-in-Aid for Scientific Research (S) (No. 22H04942) and fund (No.~215001) from the Central Research Institute of Fukuoka University.

\end{document}